\documentclass[12pt,letterpaper]{article}
\usepackage{graphicx}
\usepackage{amsmath}
\usepackage{amsfonts}
\usepackage{amsthm}
\usepackage{amssymb}
\usepackage{color}

\usepackage{verbatim} 
\usepackage{hyperref}
\usepackage{breakurl}

\newtheorem{prop}{Proposition}

\newtheorem{thm}[prop]{Theorem}
\newtheorem{lemma}[prop]{Lemma}

\theoremstyle{remark}

\newcommand{\be}{\begin{equation}}
\newcommand{\ee}{\end{equation}}
\newcommand{\bea}{\begin{eqnarray}}
\newcommand{\eea}{\end{eqnarray}}
\newcommand{\beas}{\begin{eqnarray*}}
\newcommand{\eeas}{\end{eqnarray*}}

\begin{document}

\title{Baron M\"{u}nchhausen Redeems Himself: Bounds for a Coin-Weighing Puzzle}

\author{Tanya Khovanova\\MIT \and Joel Brewster Lewis\\MIT}

\maketitle

\begin{abstract}
We investigate a coin-weighing puzzle that appeared in the Moscow Math Olympiad in 1991.  We generalize the puzzle by varying the number of participating coins, and deduce an upper bound on the number of weighings needed to solve the puzzle that is noticeably better than the trivial upper bound. In particular, we show that logarithmically-many weighings on a balance suffice.
\end{abstract}

\section{Introduction}

\subsection{Background}

Baron M\"{u}nchhausen is famous for telling the truth, only the truth and nothing but the truth \cite{Raspe}. Unfortunately, no one believes him. Alexander Shapovalov gave him an unusual chance to redeem himself by inventing a problem that appeared in the Regional round of the All-Russian Math Olympiad in 2000 \cite{BaronOriginal}. 

\begin{quote}
Eight coins weighing $1, 2,\ldots, 8$ grams are given, but which
weighs how much is unknown.  Baron M\"{u}nchhausen claims he knows
which coin is which; and offers to prove himself right by conducting one weighing on a balance scale, so as to unequivocally demonstrate the weight of at least one of the coins.  Is this possible, or is he exaggerating?
\end{quote}

In \cite{Baron}, T. Khovanova, K. Knop and A. Radul considered a natural generalization of this problem.  They defined the following sequence, which they called \emph{Baron M\"{u}nchhausen's sequence}:

\begin{quote}
Let $n$ coins weighing $1, 2,\ldots, n$ grams be given.  Suppose Baron M\"{u}nchhausen knows which coin weighs how much, but his audience does not.  Then $b(n)$ is the minimum number of weighings the Baron must conduct on a balance scale, so as to unequivocally demonstrate the weight of at least one of the coins.
\end{quote}
They completely described the sequence. Namely, they proved that $b(n) \leq 2$, and provided the list of $n$ for which $b(n) = 1$.

A similar coin-weighing puzzle, due to Sergey Tokarev \cite{tokarev}, appeared in the last round of the Moscow Math Olympiad in 1991:

\begin{quote}
You have 6 coins weighing 1, 2, 3, 4, 5 and 6 grams that look the same, except for their labels. The number $\{1, 2, 3, 4, 5, 6\}$ on the top of each coin should correspond to its weight. How can you  determine whether all the numbers are correct, using the balance scale only twice?
\end{quote}

Most people are surprised to discover that only in two weighings the weight of the each coin can be established. We invite the reader to try this puzzle out before the enjoyment is spoiled on page~\pageref{puzzle-solution}.

\subsection{The Sequence}

We generalize the preceding puzzle to $n$ coins that weigh $1$, $2$, $\ldots$, $n$ grams.  We are interested in the minimum number of weighings $a(n)$ on a balance scale that are needed in order to convince the audience about the weight of all coins.  

In this paper, we demonstrate that we can do this in not more than order of $\log n$ weighings. Because the sequence $a(n)$ relates to the task of identifying \emph{all} coins (while the sequence $b(n)$ relates to the task of identifying \emph{some} coin) we will call it \emph{the Baron's omni-sequence}.  We also calculate bounds for how many weighings are needed to prove the weight for a given particular coin.

\subsection{The Roadmap}

In Section~\ref{sec:sequence} we give a precise definition of the Baron's omni-sequence and calculate its first few terms.  In Section~\ref{sec:naturalbounds} we prove natural lower and upper bounds for the sequence, and in Section~\ref{sec:moreterms} we present the values of all known terms of the sequence. Section~\ref{sec:notation} is devoted to useful notations and terminology. 

In Section~\ref{sec:idea} we describe the idea behind the main proof of a tighter upper bound. We put this idea into practice in the subsequent three sections: in Section~\ref{sec:helpercoins}, we show that it is possible to determine the weights of several special coins in $\lceil \log_2 n\rceil$ weighings, and in Section~\ref{sec:theorem} we show how to use $\lceil \log_2 n \rceil$ aditional weighings to prove the weights of the rest of the coins.  We thus establish that $a(n)$ does not exceed $2\lceil \log_2 n \rceil$.  In Section~\ref{sec:refined}, we give a refined version of the argument which results in a modestly improved bound. 

In Section~\ref{sec:particularcoins} we consider the related task of proving the weight of a particular (e.g., adversarially-chosen) coin and prove that it can be done in not more than seven weighings.

In Section~\ref{sec:discussion} we discuss three topics. First, we discuss the question of the monotonicity of the Baron's omni-sequence. We do not come to a conclusion, but just provide considerations. Second, we show how Konstantin Knop and his collaborators used the rearrangement inequality to find optimal sets of weighings for a number of different values of $n$. Third, we give an idea of how the lower bound might be improved.

Finally, in Section~\ref{sec:futureresearch} we offer some further comments, questions and ideas for future research.

\section{The Sequence}\label{sec:sequence}

The sequence $a(n)$ is defined as follows:

\begin{quote}
Let $n$ coins weighing $1,\ 2,\ \ldots,\ n$ grams be given.  Suppose Baron M\"{u}nchhausen knows which coin weighs how much, but the audience
does not.  Then $a(n)$ is the minimum number of weighings he must conduct on a balance scale, so as to unequivocally demonstrate the weight of all the coins.
\end{quote}
The original Olympiad puzzle is asking for a proof that $a(6) = 2$.

\subsection{Examples}\label{sec:examples}

Let us see what happens for small indices. 

For $n=1$, the Baron does not need to prove anything, as there is just one coin weighing 1 gram.

For $n=2$, one weighing is enough.  The Baron places one coin on the left pan of the scale and one on the right, after which everybody knows that the lighter coin weighs 1 gram and the heavier coin weighs 2 grams.

For $n=3$, by exhaustive search we can see that the Baron can not prove all the coins in one weighing, but can in two.  For the first weighing, he compares the 1-gram and 2-gram coins, and for the second weighing the 2-gram and the 3-gram coins. Thus he establishes the order of the weights.

For $n=4$, two weighings are enough.  First, the Baron places the 1-gram and 2-gram coins on the left pan and the 4-gram coin on the right pan.  The only way for one coin to be strictly heavier than the combination of two others is for it to be the 4-gram coin.  The 3-gram is also uniquely identified by the method of elimination. In the second weighing, the Baron differentiates the 1-gram and the 2-gram coins.

For $n=5$, two weighings are enough.  The Baron places the 1-gram and 2-gram coins on the left pan and the 4-gram coin on the right pan. For the second weighing he places the 1-gram and the 4-gram coins on the left pan and the 5-gram coin on the right pan. It is left to the reader to check that these two weighings identify each coin.

For $n=6$,\label{puzzle-solution} two weighings are enough. The first weighing is $1+2+3 = 6$. This identifies the 6-gram coin and divides the other coins into two groups: $\{1,\ 2,\ 3\}$ and $\{4,\ 5\}$. The second weighing is $1+6 < 3+5$.

Another essentially different solution for $n=6$ was suggested by Max Alekseyev in a private email: $1+2+5 < 3 + 6$ and $1 + 3 < 5$.

So the sequence begins, 0, 1, 2, 2, 2, 2.

Because it is something of a mouthful to always refer to the good Baron M\"{u}nchhausen, we suppress further mention of him.  Instead, ``we'' will perform the weighings, or they will take place in the passive voice.

\section{Natural Bounds}\label{sec:naturalbounds}

For all $n$, we have that $a(n) \le n-1$ (see \cite{blog1}): for each $k < n$, in the $k$-th weighing we compare the $k$-gram and $(k+1)$-gram coins.  Getting the expected result every time confirms the weights of all coins.

On the other hand, we have that $a(n) \ge \log_3(n)$. Indeed, suppose we conduct several weighings; then to every coin we can assign a sequence of three letters $L$, $R$, $O$, corresponding to where the coin was placed during each weighing -- on the left pan, on the right pan or in the out-pile (i.e., on neither pan). If two coins are assigned the same letters for every weighing, then our weighings do not distinguish between them.  That is, if we switched the weights of these two coins, the results of all the weighings will be the same.  If the number of weighings were less than $\log_3(n)$, we are guaranteed to have such a pair of coins.  Thus, at least $\log_3(n)$ weighings are needed.

\section{More Terms}\label{sec:moreterms}

Several other terms of the sequence are known. In the cases $n = 10$ and $n = 11$, Alexey Radul found sets of three weighings that demonstrate the identity of every coin \cite{blog1}. As this matches the lower bound, we conclude that $a(10) = a(11) = 3$.  Max Alekseyev wrote a program to exhaustively search through all possible combinations of weighings, with the result that $a(7) = a(8) = a(9) = 3$. The program also confirmed the values for $n = 10$ and $n = 11$, but larger values of $n$ were beyond its limits.

After that Konstantin Knop calculated more terms of the sequence by finding weighings that match the lower bound. In particular, he stated that he found weighings to demonstrate that  $a(12) = \ldots = a(17) = 3$ and $a(53) = 4$ (see comments in \cite{blog1}). When we were writing this paper we asked Konstantin Knop if he would share his weighings with us. He sent them to us, explaining that they were calculated together with Ilya Bogdanov.  With his permission we include some of his weighings in this paper. 

Here we show how to demonstrate the identities of all coins for $n=15$. The technique is similar to the one used in cases for $n=4$ and $n=6$, and we will use a related technique in Section~\ref{sec:theorem} to prove our upper bound. 

The first weighing is
$$
1 + \dots + 7 < 14 + 15.
$$
The only way a collection of seven coins can be lighter than two coins is if the seven coins are the lightest coins from the set and the two coins are the heaviest. Thus, this weighing divides all coins into three groups $C_1 = \{1, 2, 3, 4, 5, 6, 7\}$, $C_2 = \{8, 9, 10, 11, 12, 13\}$ and $C_3 = \{14, 15\}$.

In the second weighing, the audience sees three coins from $C_1$, one coin from $C_2$ and both coins from $C_3$ go on the left pan, while three coins from $C_1$ and two coins from $C_2$ go on the right pan:
$$
(1  +  2  +  3)  + 8 + (14 + 15) = (5 + 6 + 7) + (12 + 13).
$$
Observing that the weighing balances, the audience is forced to conclude that the left pan holds the lightest coins from each group and the right pan holds the heaviest. Thus, the coins are split into the following groups: $\{1, 2, 3\}$, $\{4\}$, $\{5, 6, 7\}$, $\{8\}$, $\{9, 10, 11\}$, $\{12, 13\}$ and $\{14, 15\}$.

Similarly, we take the third weighing
$$
1  + 5 + 8 + 9 + 12 + 14 = 3 + 7 + 11 + 13 + 15,
$$
and this can balance only if the the lightest coins from each group are on the left pan and the heaviest are on the right.  Thus, in the end all coins are identified.

The other weighings that Konstantin Knop sent to us use a different technique which is not related to our proof of the upper bound for $a(n)$, so we delay presenting it until Section~\ref{sec:rearrangement}.  Maxim Kalenkov used the same technique and the help of a computer to find two more terms, namely $a(18) = a(19) = 3$.

So the sequence begins, 0, 1, 2, 2, 2, 2, 3, 3, 3, 3, 3, 3, 3, 3, 3, 3, 3, 3, 3.

No sets of three weighings that identify all coins are known for $20 \leq n \leq 27$.  However, Maxim Kalenkov found solutions in four weighings for a range of numbers from $n=20$ to $n=58$ inclusive.

\section{Notation and Terminology}\label{sec:notation}

For integers $x \leq y$, we denote by $[x \dots y]$ the set of consecutive integers between $x$ and $y$, inclusive.  
For $x=1$, instead of using $[x\dots y]$ we will just use $[y]$, which is the standard notation for the range anyway.

We will use the number $i$ to denote the $i$-gram coin on a pan.  Thus, $[x \dots y]$ represents the set of coins 
of weights no smaller than $x$ and no larger than $y$, and we will occasionally construct weighings using this set 
notation.  All arithmetic operations \emph{other than addition} are understood to take place on the weight of a single 
coin; thus $3 \cdot 2^2 - 1$ represents the $11$-gram coin.  Inside brackets, addition operates in the same way, 
so $[3+4\ldots 11-1]$ is the set of coins weighing from $7$ to $10$, inclusive, and does not include the coins $3$, 
$4$, $11$, or $1$. Outside square brackets, addition can be read to have the same meaning as union, so $1 + 2$ means 
the $1$-gram and the $2$-gram coins taken together, while $[3] + [5 \dots 7]$ is the set $\{1, 2, 3, 5, 6, 7\}$.

Equalities and inequalities represent the outcomes of particular weighings; thus $[3] + [5 \dots 7] > [11 \dots 12]$ represents the weighing with the coins 1, 2, 3, 5, 6, and 7 on the left pan and the coins 11 and 12 on the right pan, in which the left pan had larger total weight.  In particular, when representing a weighing as an equality/inequality we will refer to the left and right sides of the equality/inequality as the left and
right pans of the weighing, respectively.

If $A$ denotes a set of coins, then $|A|$ denotes the total weight of
those coins.  Note that this is \emph{not} the cardinality of the set $A$, which we denote $\# A$.

Define the \emph{small half} of a (finite, totally ordered) set $A$ to be the set consisting of the $\left\lfloor \frac{\# A}{2} \right\rfloor$ smallest elements of $A$.  For example, the small half of $\{1, 3, 4, 5, 7\}$ is $\{1, 3\}$.

A subset $B$ of a set $A$ is said to be \emph{upwards-closed} if for every $x \in B$ and $y \in A$ with $x < y$ we have $y \in B$.  Thus, the set $\{1, 3, 4, 5, 7\}$ has six upwards-closed subsets, three of which are $\{4, 5, 7\}$, the entire set, and the empty set.  The notion of a \emph{downwards-closed} subset is defined analogously.

\section{An Idea}\label{sec:idea}

Before we proceed with the main section of the proof, we present an idea that actually does not work, but that we will use as a starting point.

Given a set of coins $[1\ldots n]$, suppose we can find numbers $k < m$ such that $|[1\ldots k]| = |[m+1\ldots n]|$. In this case the weighing
\[
[1\ldots k] = [m+1\ldots n]
\]
will balance. This fact demonstrates that the coins in question really are the coins we claim: the sum of $k$ coins is at least the weight of the left pan, while the sum of $n - m$ coins is at most the weight of the right pan.  This gives us our first division into three parts: $[1\ldots k]$, $[k+1\ldots m]$ and $[m+1\ldots n]$. 

If we have in particular that $k = n/2$, then the division into the three parts above supplies us with the division of the range $[n]$ into two halves.

Suppose for the second weighing we can balance $[1\ldots n/4] + [n/2+1\ldots 3n/4]$ against an appropriately-chosen combination of upwards-closed subsets of $[n/4 + 1\ldots n/2]$ and $[n/2 + 1\ldots n]$.  In this way we divide each half from the previous division into halves again. 

For the third division, we place the small half of each of the four groups into which we have divided the coins on one pan, and we choose an upwards-closed subset of the heavier halves on the other pan so that the pans balance.  This again divides each of our four subsets into two halves.

Continuing such binary division we can identify all coins in $\log_2(n)$ weighings.

Unfortunately, this strategy fails in a very simple way: it is impossible to carry out in general.  In particular, the very first step is quite often impossible. Consider, for example, 12 coins. We want to find an upward-closed subset to balance the lightest six coins. But $12 < 1 + 2 + 3 + 4 + 5 + 6 < 12 + 11$. 

However, this problem can be overcome if we are first able to prove the identities of a small number of helper coins; these coins could then be used to make up the difference between the small half of the coins and the corresponding upward-closed set.

For example, if we start with 12 coins and somehow can prove the identities of the 1-gram and the 2-gram coins, then we can balance out the small half of the leftover set: $3+4+5+6+7 = 12 + 11 + 2$.  This suggests that we should start by looking for easily-identifiable sets of ``helper coins.''

\section{Helper Coins}\label{sec:helpercoins}

We now give a simple procedure to identify a set of helper coins.  This set of helper coins does not require many weighings to identify.  In addition, it is versatile and produces many sums.

Let the binary expansion of $n - 1$ be $n - 1 = 2^{a_1} + 2^{a_2} + \ldots $ with $a_1 > a_2 > \ldots \geq 0$ and $a_1 = \lceil \log_2 n \rceil - 1$.  We perform the weighings $1 < 2$, $1 + 2 < 4$, $1 + 2 + 4 < 8$, \ldots, $1 + 2 + 4 + \ldots + 2^{a_1 - 1} < 2^{a_1}$ and $2^{a_1} + 2^{a_2} + \ldots < n$.  From the first weighing, we learn that the coin we claim has weight $2$ grams has weight at least this large.  Similarly, from the second weighing we learn that the coin we claim has weight $4$ grams weighs at least that much, and so on.  Thus, the last weighing demonstrates that the coin we claim has weight $n$ has weight at least $n$.  However, all our coins weigh at most $n$ grams, whence the coin we claim has weight $n$ must actually have that weight.  Moreover, this also shows that the coins $1$, $2$, \ldots, $2^{a_1}$ are the coins we claim.

Denote by $S(n)$ the set of coins identified by these $\lceil \log_2 n \rceil$ weighings.  The following useful property of $S(n)$ is clear.
\begin{prop}\label{prop:usefulproperty}
Using only the elements of $S(n)$, we can construct a pile of coins whose weight is $i$ for any $i \in [n]$.  That is, $[n] \subseteq \{|T| \colon T \subseteq S(n)\}$.
\end{prop}

We now use this set $S(n)$ to give an effective version of the algorithm described in Section~\ref{sec:idea}.

\section{The Upper Bound}\label{sec:theorem}

\begin{thm}
We can identify all coins in $[n]$ in at most $2\lceil \log_2 n \rceil$ weighings.  That is, $a(n) \leq 2\lceil \log_2 n \rceil$.
\end{thm}

\begin{proof}
Section \ref{sec:examples} demonstrates the result for $n = 1, 2, 3$.  For $n \geq 4$, use the construction of Section \ref{sec:helpercoins} to identify the coins in the set $S(n)$ in $\lceil \log_2 n\rceil$ weighings. 

Set $C = [n]\smallsetminus S(n)$.  We perform binary search on $C$ as follows: suppose that at some stage, we have successfully demonstrated a division of $C = C_{1} \cup C_{2} \cup \cdots \cup C_{m}$ into several disjoint ranges so for every non-helper coin we know to which range it belongs, and that the ranges are numbered in order: for any $i < j$ and any $x \in C_i,\ y \in C_j$ we have $x < y$. (Initially, this is the case with $m = 1$ and $C = C_1$.)

If $C_i$ consists of one element, then the identity of the coin in $C_i$ is already proven, so we may set it aside.  For all $i$ for which $\# C_i > 1$, we split $C_i$ and place its small half on the left pan of the balance.  The other non-proven elements of $C$ have larger total weight and each is of weight at most $n$.  Thus, we may begin adding unused non-proven elements of $C$ to the right pan, starting with the largest, until the right pan weighs at least as much as the left pan. 
As soon as the right pan reaches the weight of the left pan the difference between the weights of the two pans is at most $n$.  Then (by Proposition \ref{prop:usefulproperty}) we may add elements from $S(n)$ to the left pan as needed in order to make the pans balance.  

This weighing identifies the small half of $C_i$ for all $i$, and so divides each $C_i$ into two almost-equal-sized parts.  Repeating this $\log_2 (\# C) \leq \lceil \log_2 n \rceil$ times results in a total ordering of the elements of $C$, and so the identification of all elements of $[n]$.  Thus, at most $2 \lceil \log_2 n \rceil$ weighings are required.
\end{proof}

\section{The Refined Upper Bound}\label{sec:refined}

The way we divide coins into piles in the previous theorem is not optimal.  In particular, it leaves room for two improvements.  First, when we remove the small half of each pile (to place on the left pan), we ignore the information we get from the fact that the remaining coins are divided into two parts (some on the right pan, some in the out pile).  Thus, we can do better by keeping track of all three parts of the division.  Second, the most profitable way of dividing coins into three piles would be to have each pile of the same, or almost the same, size. In our approach it is not possible for the set of the lightest coins to be the same cardinality as the set of the heaviest coins of the same total weight.  However, it is possible to do better than in Section~\ref{sec:theorem} by choosing a division in which the part of largest cardinality has less than half of the coins.

Suppose we have the set of coins $[n]$.  For some $k, m$, we divide the coins between the two pans (with some left out, i.e., not on either pan) by placing the lightest $k$ coins on the left pan and the heaviest $m$ coins on the right pan so that the right pan holds more total weight and the weight difference between the two pans does not exceed $n$.  In this case all coins are divided into three groups of sizes $k$, $m$ and $n-k-m$.  We seek values of $k, m$ that give an optimal division of this form.

\begin{lemma}
Subject to our conditions, in an optimal division we have that no pile contains more than $(-2 + \sqrt{6})n$ coins.
\end{lemma}

\begin{proof}
The lightest $k$ coins weigh slightly more than $k^2/2$ grams, so the pile on the right pan weighs at least $k^2/2$ grams.  As each coin weighs not more than $n$ grams, it follows that the pile on the right pan has at least $k^2/(2n)$ coins.  Hence the out-pile contains not more than $n-k-k^2/(2n)$ coins. As the right pan is guaranteed to have fewer coins than the left pan, to build an optimal division we need to have the same number of coins on the left pan as out, and thus the optimal choice of $k$ satisfies $k  \leq n-k-k^2/(2n)$.  Consequently, the value of $k$ that satisfies $k = n-k-k^2/(2n)$ will be no smaller than (and presumably close to) the optimal value.  Elementary algebra gives the result.
\end{proof}

Define $\alpha = -2 + \sqrt{6}$.  Observe that whether $k = \alpha n$ is actually optimal or just slightly larger than optimal, we have that this value of $k$ is better than choosing $k = n/2$ (as in Section~\ref{sec:theorem}). 

Recall that in practice, we divide into piles not the full range $[n]$, but rather a set $C$ from which our helper coins are excluded. In the following lemma we prove that we can keep our new estimate for such a set.

\begin{lemma}
Given a subset $C'$ of $[n]$, there exists a partition of $C'$ into three parts meeting the following conditions: the first part consists of some lightest coins in $C'$ (i.e., it is downwards-closed) and the second consists of some of the heaviest; the subset with heaviest coins weighs more than the one with the lightest, but not by more than $n$; and none of the three parts contains more than $\alpha n$ coins.
\end{lemma}

\begin{proof}
As a first approximation of the desired weighing, we place the $\alpha n$
smallest coins from $C'$ on the left pan.  On the right pan, we place the smallest possible upwards-closed subset of $C'$ such that the right pan is not lighter than the left pan and the difference of their weights is not more than $n$.  In this case, it is clear that the right pan can not have more coins than the left pan. In addition, we know that the left pan has more total weight than the weight of the smallest $\alpha n$ coins from the range $[n]$, and we formerly required the $(1 - 2\alpha)n$ heaviest coins from the range $[n]$ to overbalance the left pan; since some of the heaviest coins in $[n]$ might be missing from $C'$, we might need even more than $(1 - 2\alpha)n$ coins in the right pan.  Hence, the right pan has at least $(1 - 2\alpha)n$ coins and so the out pile will have not more than $\alpha n$ coins.

The only problem we can encounter is that we can run out of coins for the right pan before the right pan reaches the weight of the left pan. In this case we perform the following procedure. We remove the heaviest coin from the left pan. If this is enough to have the balance we need, we are done. If this is not enough, we place that removed coin on the right pan. We continue until we get a weighing satisfying our requirements. At the end of this process, the left pan can not
have more than $\alpha n$ coins and the right pan can not have more coins than
the left pan, and the out pile in this case will be not more than one coin.
\end{proof}

To finish what we have started, we need to remember that not only the first division into piles needs to be optimal. We continue with subdivisions. Intuitively, in every subsequent step it should be easier to form balanced divisions, because the coins in each subset have a smaller spread of weights. The lemma below guarantees that we can continue the divisions in such a way that the maximum pile size at every next step will not exceed $\alpha$ times the maximum pile size at the previous step.

\begin{lemma}\label{lemma:subdivisions}
Given a set of coins whose weights are distinct positive integers between $a$ and $b$, we can divide it into three groups, the lightest, the heaviest and the middle, so that the following holds: the size of each group is not more than $\lceil \alpha (b-a+1) \rceil$, the second group weighs more than the first group, and the difference between the weights of these two groups is not more than $b-a+1$.
\end{lemma}
\begin{proof}
The proof is the same as the previous proof, \emph{mutatis mutandis}.
\end{proof}

This refined approach gives us a better upper bound.

\begin{thm}
We can identify all coins in $[n]$ in at most $\lceil \log_2 n\rceil + \lceil \log_{\alpha^{-1}} n \rceil$ weighings.
\end{thm}

For comparison, the bound of Section~\ref{sec:theorem} is approximately $2 \log_2 n \approx 3.17 \log_3 n$ while our new bound is about $2.96 \log_3 n$.

\section{Particular Coins}\label{sec:particularcoins}

One of the future research questions in \cite{Baron} was to find the minimum number of weighings needed if the audience requests that the Baron prove the weight of a particular coin.  For our purposes, it is tempting to think that for all $n$ (or at least for sufficiently large $n$), some particular coin $t(n)$ might require order of $\log n$ weighings to identify, and so perhaps give an improvement over the lower bound of Section \ref{sec:naturalbounds}.  The following theorem rules out this possibility.  In particular, we show that for each positive integer $t$ and for any $n \geq t$, the coin of weight $t$ can be identified among the coins $[n]$ in at most seven weighings.  Our proof relies on the following number-theoretic property of triangular numbers proved by Gauss \cite{Gauss, Grosswald}.

\begin{lemma}\label{lemma:triangular}
Every positive integer $n$ can be written as the sum of three (not necessarily distinct) triangular numbers, possibly including $0$.
\end{lemma}

For notational convenience, we denote by $T_\ell$ the $\ell$-th triangular number $T_{\ell} = \frac{\ell(\ell + 1)}{2}$.

\begin{thm}
Given any $t \in [n]$, we can identify the coin $t$ in seven weighings.
\end{thm}

\begin{proof}
The result is true for small values of $n$ either by the results of Section~\ref{sec:moreterms} or from our upper bounds on the Baron's omni-sequence, so suppose $n > 8$.

First, we show that for most values of $t$ we can identify the $t$-coin in only six weighings.  In particular, suppose that $t \geq \sqrt{2n}$.

By Lemma \ref{lemma:triangular}, there exist integers $a \leq b \leq c$ such that $t = T_a + T_b + T_c$. If $a > 0$, then we perform the three weighings
\[
[1\ldots c] = T_c,
\]
\[
[1\ldots b] + T_c = (t - T_a),
\]
and
\[
[1\ldots a] + (t - T_a) = t
\]
each with exactly one coin on the right pan of the balance.  From these weighings, we may conclude that the coins that appear on the right pan weigh at least as much as we claim, and in particular that the coin $t$ weighs at least $t$ grams.  If $a = 0$ or $a = b = 0$, then we omit respectively the third weighing or the second and third weighings, and have the same conclusion. 

Similarly, there exist integers $i \leq j \leq k$ such that $n - t = T_i + T_j + T_k$.  If $i > 0$, then we perform the three weighings
\[
[1\ldots k] + t = (T_k + t),
\]
\[
[1\ldots j] + (T_k + t) = (n - T_i),
\]
and
\[
[1\ldots i] + (n - T_i) = n,
\]
concluding that the only way all these weighings can balance is for each weighing to have on the right pan the coin that we claim is on the right pan. In particular, the coin $t$ has to weigh exactly $t$ grams, as needed. As above, if $i = 0$ then we may reach the same conclusion using fewer weighings.

Observe that $T_k < n$ and so $k < \sqrt{2n} \leq t$, so all of the weighings above use each coin at most once.  However, if we choose $t < \sqrt{2n}$ then in the fourth weighing we might need two coins that weigh $t$ grams: once as one coin in the range $[k]$ and also as an individual coin.  Thus, we need some other scheme of weighings in this case.  For these lighter choices of $t$ we replace the last three weighing with the four weighings
\[
[1\ldots k] = T_k,
\]
\[
[1\ldots j] + T_k = (n - t - T_i),
\]
\[
[1\ldots i] + (n - t - T_i) = (n - t),
\]
and
\[
(n - t) + t = n.
\]
From these weighings, we may conclude that the coins that appear on the right pan weigh at least as much as we claim, and in particular that the coin $n$ weighs at least $n$ grams.  However, no coin weighs more than $n$ grams, so actually the inequalties deduced from the first six weighings must actually be equalities.  In particular, the coin we claimed to be of weight $t$ does weigh exactly $t$ grams, as needed. 

(We note that the supposition $n > 8$ is necessary only to insure that $t < \sqrt{2n}$ implies $t < n - t$, so the coins $t$ and $n - t$ are different and the final weighing is legal.)
\end{proof}

\section{Discussion}\label{sec:discussion}

In this section we discuss a variety of topics related to the Baron's omni-sequence, including the question of monotonicity of the sequence, methods of finding optimal sets of weighings, and how to improve the lower bound.

\subsection{Is the Sequence Non-Decreasing?}\label{sec:non-decreasing}

There is no reason to believe that the sequence $a(n)$ is non-decreasing. The possibility of dividing coins into nice weighings often appears to depend on properties of the particular number $n$. It is conceivable that from time to time we encounter a number $n$ such that it is easier (i.e., requires fewer weighings) to identify all the coins of $[n]$ than all the coins of $[n - 1]$.

In addition, there are other sequences related to the sequence $a(n)$ that are not monotonic. For example, the sequence $b(n)$ of the minimum number of weighings Baron needs to prove one coin of his choosing is completely described in \cite{Baron} and it fluctuates between the values 1 and 2.

For $t \leq n$, denote by $m_t(n)$ the minimum number of weighings necessary to identify the coin $t$ among the set of coins $[n]$.  In Section~\ref{sec:particularcoins}, we showed that $m_t(n) \leq 7$ for all $n, t$.  Here, we present another statement regarding the behavior of this function.  

\begin{prop}
There exists a number $t$ such that the sequence $m_t(n)$ (of the minimum number of weighings necessary to prove the identity of the coin $t$) is not monotonic.
\end{prop}
\begin{proof}
From the paper \cite{Baron} we know that the values of $n$ such that there exists a coin that can be found in one weighing belong to one of the following four groups: 
\begin{itemize}
\item If $n$ is a triangular number then the $n$-coin can be found in one weighing.  These values of $n$ are sequence A000217 in the OEIS \cite{OEIS}: 1, 3, 6, 10, 15, 21, \ldots.
\item If $n$ is one more than a triangular numbers then the $n$-coin can be found in one weighing.  These values of $n$ are sequence A000124 in the OEIS \cite{OEIS}: 2, 4, 7, 11, 16, 22, \ldots.
\item If the $n$-th triangular number is one more than a perfect square, i.e., if $T_n=k^2+1$ for some integer $k$, then the $k$-coin can be found in one weighing.  These values of $k$ are sequence A106328 in the OEIS \cite{OEIS}: 3, 18, 105, 612, 3567, 20790, \ldots.
\item If the $n$-th triangular number is a perfect square, i.e., if $T_n=k^2$ for some integer $k$, then the $k$-coin can be found in one weighing.  These values of $k$ are sequence A001109 in the OEIS \cite{OEIS}: 1, 6, 35, 204, \ldots.
\end{itemize}

Now we see that in particular, $m_{105}(n)$ is not monotonic.  Since $105$ is a triangular number we have that $m_{105}(105) = 1$.  Also, since $T_{148} = 11026 = 105^2+1$ we have that $m_{105}(148) = 1$.  However, for any $n$ between $105$ and $148$ we have $m_{105}(n) \geq 2$.  Thus the sequence $m_{105}(n)$ is not monotonic.
\end{proof}

\subsection{The Rearrangement Inequality}\label{sec:rearrangement}

All the solutions in the cases $n > 15$ that were sent to us by Konstantin Knop were proved using a technique that we have not yet mentioned in this paper.  In particular, they make use of the following classical inequality:
\begin{lemma}[Rearrangement inequality]
Given two sets of distinct real numbers $a_1 < a_2 < \ldots < a_n$ and $b_1 < b_2 < \ldots < b_n$.  As $\sigma$ varies over the permutations of $[n]$, the value 
$a_1 b_{\sigma(1)} + a_2 b_{\sigma(2)} + \ldots + a_n b_{\sigma(n)}$ achieves its minimum for the reverse-identity permutation, i.e., when $\sigma(i) = n-i+1$ for all $i$.  Furthermore, this minimum is unique.
\end{lemma}

The following weighings were found by Maxim Kalenkov. Let $n =19$, and for clarity let $c_i$ denote the $i$-gram coin.
Consider the weighings
$$c_1+c_2+c_3+c_4+c_5+c_7+c_8+c_{10}+c_{13} = c_{16}+c_{18}+c_{19},$$
$$c_1+c_2+c_3+c_6+c_9+c_{11}+c_{16} = c_8+c_{10}+c_{13}+c_{17},$$
and
$$c_1+c_4+c_6+c_8+c_{12}+c_{18} = c_3+c_7+c_{11}+c_{13}+c_{15}.$$
Multiply the first equation by $12$, the second by $7$, the third by $3$, and sum them up.  We get the following equality: 
\begin{multline*}
22c_1 + 19c_2 + 16c_3 + 15c_4 + 12c_5 + 10c_6 + 9c_7 +
8c_8 + 7c_9 + 5c_{10} + \\
 4c_{11} + 3c_{12} + 2c_{13} + 0c_{14} -3c_{15} -5c_{16} -7c_{17} -9c_{18} -12c_{19} = 0.
\end{multline*}

By the rearrangement inequality, the lowest value the left-hand side of this equation can achieve is when coins match their labels and this value is $22 \cdot 1 + 19 \cdot 2 + 16 \cdot 3 + 15 \cdot 4 + 12 \cdot 5 + 10 \cdot 6 + 9 \cdot 7 +
8 \cdot 8 + 7 \cdot 9 + 5 \cdot 10 + 4 \cdot 11 + 3 \cdot 12 + 2 \cdot 13 + 0 \cdot 14
-3 \cdot 15 -5 \cdot 16 -7 \cdot 17 -9 \cdot 18 -12 \cdot 19$ which is equal to zero. Thus, the only way to achieve equality is for every coin to match the labeling.

A similar approach was used to find a solution in three weighings for $n = 16$, $17$, $18$ and $19$, as well as solutions in four weighings for $20 \leq n \leq 58$.

\subsection{The Lower Bound}\label{sec:lowerbound}

So far, we have not discussed improvements on the natural lower bound, mostly because all the examples we know are very close to it.  Since currently the only method we have to find the terms of the Baron's omni-sequence is to find weighings that match the lower bound, even a slight improvement in the lower bound can be extremely useful.  Here we show an idea of how we might be able to improve the lower bound.

For every $n \leq 19$, there is an optimal set of weighings such that one weighing is very special.  First, in this weighing either the two pans balance or the lighter pan is exactly one gram lighter than the heavier pan.  Second, in this weighing all the coins on the left pan are lighter than all the coins on the right pan.

Suppose that one could prove that for every $n$, there exists an optimal set of weighings such that one of the weighings has the two properties above.  Under this assumption, this particular weighing can not divide the coins into three equally-sized groups. In particular, it would follow that the number of weighings for $n=9$ has to be more than two, and for $n=27$, more than three.

Let us estimate the bounds on the sizes of the three groups of coins that can be achieved during this special weighing. Suppose the left pan has $\gamma n$ coins and the right pan has $\delta n$ coins, with every coin on the left lighter than every coin on the right. The $k$-th lightest coin on the right pan weighs at least $\gamma n$ more grams than the $k$-th lightest coin on the left pan, so the coins on the right pan weigh $\gamma \delta n^2$ more than the lightest $\delta n$ coins from the left pan. Each leftover coin on the left pan weighs not more than $(1 - \delta)n$ grams. Hence, the left pan has at least $\frac{\gamma \delta}{1 - \delta} n$ more coins than the right pan and so we have 
\begin{equation}\label{eq:lowerbound}
\gamma \geq \delta + \frac{\gamma \delta}{1 - \delta}.
\end{equation}
Thus, we seek to minimize the quantity $\max(\gamma, \delta, 1 - \gamma - \delta)$ subject to Equation \ref{eq:lowerbound} and the constraints that $\gamma$, $\delta$ and $1 - \gamma - \delta$ are nonnegative.  Under these conditions, this best possible division into three piles is achieved when 
\[ 
\gamma = \delta + \frac{\gamma \delta}{1 - \delta}
\]
and
\[ 
\gamma = 1 - \gamma - \delta,
\]
at which point we have $\max(\gamma, \delta, 1 - \gamma - \delta) = \gamma = \frac{3}{8}$.  So, conditional upon the existence of such a special weighing in an optimal set of weighings, we have that the lower bound may be raised to $\lceil \log_3(3n/8) \rceil + 1$.  This new estimate implies that we will need four weighings for $n \geq 25$ and five weighings for $n \geq 73$.

\section{Future Research}\label{sec:futureresearch}

Our upper bound can be improved by tightening Lemma \ref{lemma:subdivisions}.  For example, it should be possible to show that starting from the second step our dividing constant can be better than $\alpha$.  The best-case for such an argument (based on the helper coins) is $\log_2 n + \log_3 n \approx 2.58 \log_3 n$.  

It seems like the true growth rate of the sequence may be very close to the natural lower bound of $\log_3 n$. For example, our lower bound for $n = 58$ is 4, and the refined upper bound presented in Section~\ref{sec:refined} is 12. The fact that $a(58)=4$ tells us that the lower bound estimate might be very tight.  On the other hand, if the sequence is non-monotonic then the upper and lower bounds could be quite different.  Consequently, it would be nice to determine whether or not the sequence is monotonic.

As we mentioned above, any improvement on the lower bound helps by providing hope that we can invent clever constructions to calculate more terms of the sequence.  It would be interesting to prove the conjectures we offered in Section \ref{sec:lowerbound}, namely that there always exists an optimal solution with one special weighing that fails to balance by at most one gram and in which all the coins on one pan are lighter than all the coins on the other. Is there any other way to show that the number of weighings required to identify all coins in $[n]$ is larger than the trivial lower bound? Can we prove any theorems that allow exhaustive search to become feasible for $n \geq 12$?  Or can we improve the exhaustive search and check all possibilities in a smarter way?

Using the rearrangement inequality to find good weighings seems very promising.  In the worst case, it seems like a promising way to produce sets of weighings for larger $n$ that efficiently identify all coins.  Even if the number of weighings used were to not match the lower bound exactly, it would allow for improved bounds on $a(n)$.  Is it possible to produce solutions with a small number of weighings for some infinite sequence of $n$-values?

\section{Acknowledgements}

We are grateful to Konstantin Knop and Maxim Kalenkov for sharing their weighings with us, and to all other enthusiasts who got excited by this sequence and calculated or tried to calculate some values.

%
%
%
%
%

\end{document}